\documentclass[12pt,leqno,a4paper,dvips,twoside]{article}
\usepackage{graphicx}
\usepackage{amsfonts}
\usepackage{amssymb}
\usepackage{amsmath}
\usepackage[utf8]{inputenc}
\newtheorem{theo}{Theorem}
\newtheorem{defin}{Definition}

\newtheorem{cor}{Corollary}

\begin{document}
\thispagestyle{plain} {\footnotesize {{\bf General Mathematics
Vol. xx, No. x (201x), xx--xx}}} \vspace*{2cm}

\begin{center}
{\Large {\bf A Note On Multi Poly-Euler Numbers And Bernoulli Polynomials}}
\vspace{0.5cm} \footnote{\it Received 08 Jun, 2009

\hspace{0.1cm} Accepted for publication (in revised form) 29
November, 2013}\vspace{0.3cm}

{\large Hassan Jolany, Mohsen Aliabadi, Roberto B.  Corcino, and M.R.Darafsheh }
\end{center}
\vspace{0.5cm}

\begin{abstract}
In this paper we introduce the generalization of Multi Poly-Euler polynomials and we investigate some relationship involving Multi Poly-Euler polynomials. Obtaining a closed formula for generalization of Multi Poly-Euler numbers therefore seems to be a natural and important problem. 
\end{abstract}

\begin{center}

\end{center}
\vspace*{0.3cm}


\section{Introduction}
     In the 17th century a topic of mathematical interst was finite sums of powers of
integers such as the series $1+2+...+(n-1)$ or the series $1^2 + 2^2 + ... + (n -1)^2$.The
closed form for these finite sums were known ,but the sum of the more general series
$1^k+2^k+...+(n-1)^k$was not.It was the mathematician Jacob Bernoulli who would
solve this problem.Bernoulli numbers arise in Taylor series in the expansion 
\begin{equation}\label{e1}
\begin{array}{c}
\frac{x}{e^x-1}=\sum\limits_{n=0}^{\infty}B_n\frac{x^n}{n!}
\end{array}.
\end{equation}
and we have, 
\begin{equation}\label{e1}
\begin{array}{c}S_m(n) = \sum_{k=1}^n k^m = 1^m + 2^m+ \cdots + n^m= {1\over{m+1}}\sum_{k=0}^m {m+1\choose{k}} B_k\; n^{m+1-k}
\end{array}.
\end{equation}
and we have following matrix representation for Bernoulli numbers(for $n\in \mathbf{N}$),[1-4].

\begin{align}
B_{n} &=\frac{(-1)^n} {(n-1)!}~ \begin{vmatrix}   \frac{1}{2}& \frac{1}{3} & \frac{1}{4} &\cdots \frac{1}{n}&~\frac{1}{n+1}\\
                                                             1&  1 & 1 &\cdots 1 & 1 \\
                                                             0&  2 & 3 &\cdots {n-1} & n\\
                                                               0&  0 & \binom{3}{2} &\cdots \binom{n-1}{2} & \binom{n}{2}  \\
                                                           \vdots & ~   \vdots & ~ \vdots   &\ddots~~  \vdots & \vdots &    \\
                                                             0&  0& 0& \cdots \binom{n-1}{n-2} & \binom{n}{n-2} \\
                                                             \end{vmatrix}.
\end{align}
Euler on page 499 of [5], introduced Euler polynomials, to evaluate the alternating sum
\begin{equation}\label{e1}
\begin{array}{c}
A_n(m)=\sum\limits_{k=1}^{m}(-1)^ {m-k}k^n=m^n-(m-1)^n+...+(-1)^ {m-1}1^n
\end{array}.
\end{equation}
The Euler numbers may be defined by the following generating functions
\begin{equation}\label{e1}
\begin{array}{c}
\frac{2}{e^{t}\!+\!1}\;=\;\sum\limits_{{n=0}}^{\infty}E_{n}\frac{t^{n}}{n!}
\end{array}.
\end{equation}
and we have following folowing matrix representation for Euler numbers, [1,2,3,4].

\begin{align}
E_{2n} &=(-1)^n (2n)!~ \begin{vmatrix}   \frac{1}{2!}& 1 &~& ~&~\\
                                                             \frac{1}{4!}&  \frac{1}{2!} & 1 &~&~\\
                                                                 \vdots & ~  &  \ddots~~ &\ddots~~ & ~\\
                                                               \frac{1}{(2n-2)!}& \frac{1}{(2n-4)!}& ~&\frac{1}{2!} &  1\\
                                                               \frac{1}{(2n)!}&\frac{1}{(2n-2)!}& \cdots &  \frac{1}{4!} & \frac{1}{2!}\end{vmatrix}.
\end{align}

The poly-Bernoulli polynomials have been studied by many researchers in recent decade. 
The history of these polynomials goes back to Kaneko. The poly-Bernoulli polynomials 
have wide-ranging application from number theory and combinatorics and other fields  of 
applied mathematics. One of applications of poly-Bernoulli numbers that was investigated by Chad Brewbaker in [6,7,8,9], is about the number of $(0,1)$-matrices with $n$-rows and $k $
columns. He showed the number of $(0,1)$-matrices with $n$-rows and $k$ columns uniquely 
reconstructable from their row and column sums are the poly-Bernoulli numbers of negative 
index   $B_n^{(-k)}$.
Let us briefly recall poly-Bernoulli numbers and polynomials.
For an integer $k\in \mathbf{Z}$, put  

\begin{equation}\label{e1}
\begin{array}{c}
\operatorname{Li}_k(z) = \sum_{n=1}^\infty {z^n \over n^k}.
\end{array}.
\end{equation}  
which is the $k$-th polylogarithm if $k\geq 1$ , and a rational function if $k\leq 0$. The name of the function come from the fact that it may alternatively be defined as the repeated integral of itself .
 The formal power series can be used to define Poly-Bernoulli numbers and polynomials. The polynomials $B_n^{(k)}(x)$ are said to be poly-Bernoulli polynomials if they satisfy, 
\begin{equation}\label{e1}
\begin{array}{c}
{Li_{k}(1-e^{-t}) \over 1-e^{-t}}e^{xt}=\sum\limits_{n=0}^{\infty}B_{n}^{(k)}(x){t^{n}\over n!}
 \end{array}.
\end{equation}
 In fact, Poly-Bernoulli polynomials are generalization of Bernoulli polynomials, because for $n\leq 0$, we have,

\begin{equation}\label{e1}
\begin{array}{c}
(-1)^nB_n^{(1)}(-x)=B_n(x)
\end{array}.
\end{equation} 
 Sasaki,[10], Japanese mathematician, found the Euler type version of these polynomials, In fact, he by using  the following relation for Euler numbers, 
 \begin{equation}\label{e1}
\begin{array}{c}
cosh t=\sum\limits_{n=0}^{\infty}\frac{E_n}{n!}t^n
 \end{array}.
\end{equation} 
found a  poly-Euler version as follows

\begin{equation}\label{e1}
\begin{array}{c}
\frac{Li_k(1-e^{-4t})}{4tcosht}=\sum\limits_{n=0}^{\infty}E_{n}^{(k)}{t^{n}\over n!}
\end{array}.
\end{equation} 
Moreover, he by defining the following $L$-function, interpolated his definition about Poly-Euler numbers.
\begin{equation}\label{e1}
\begin{array}{c}
L_k(s)=\frac{1}{\Gamma(s)}\int_0^{\infty}t^{s-1}\frac{Li_k(1-e^{-4t})}{4(e^t+e^{-t})}dt
 \end{array}.
\end{equation} 
and  Sasaki showed that 

\begin{equation}\label{e1}
\begin{array}{c}
L_k(-n)=(-1)^nn\frac{E_{n-1}^{(k)}}{2}
 \end{array}.
\end{equation}
But the fact is that working on such type of generating function for finding some identities is not so easy. So by inspiration of the definitions of Euler numbers and Bernoulli numbers, we can define Poly-Euler numbers and polynomials as follows which also A.Bayad [11], defined it by following method in same times.
 
\begin{defin}\label{d1}
(Poly-Euler polynomials):The Poly-Euler polynomials may be defined by using the following generating function,

\begin{equation}\label{e1}
\begin{array}{c}
\frac{2Li_k(1-e^{-t})}{1+e^t}e^{xt}=\sum\limits_{n=0}^{\infty}\mathbf{E}_{n}^{(k)}{t^{n}\over n!}
\end{array}.
\end{equation}
\end{defin}

If we replace $t$ by $4t$ and take $x=1/2$ and using the definition $cosht=\frac{e^t+e^{-t}}{2}$, we get the Poly-Euler numbers which was introduced by Sasaki and Bayad and also we can find same interpolating function for them (with some additional constant coefficient).

The generalization of poly-logarithm is defined by the following infinite series

\begin{equation}\label{e1}
\begin{array}{c}
Li_{(k_1,k_2,...,k_r)}(z)=\sum\limits_{m_1,m_2,...,m_r}\frac{z^{m_r}}{m_1^{k_1}...m_r^{k_r}}
\end{array}.
\end{equation}
which here in summation ($0<m_1<m_2<...m_r$).

Kim-Kim [12], one of student of Taekyun Kim introduced the Multi poly- Bernoulli numbers
and proved that special values of certain zeta functions at non-positive integers can be
described in terms of these numbers. The study of Multi poly-Bernoulli numbers and their combinatorial relations has received much attention in [6-13]. The Multi Poly-Bernoulli numbers may be defined as follows

\begin{equation}\label{e1}
\begin{array}{c}
\frac{Li_{(k_1,k_2,...,k_r)}(1-e^{-t})}{(1-e^{-t})^r}=\sum\limits_{n=0}^{\infty}B_n^{(k_1,k_2,...,k_r)}\frac{t^n}{n!}
\end{array}.
\end{equation}

So by inspiration of this definition we can define the Multi Poly-Euler numbers and polynomials .

\begin{defin}\label{d1} Multi Poly-Euler polynomials $\mathbf{E}_n^{(k_1,...,k_r)}(x)$, $(n = 0, 1, 2, ...)$ are defined for each integer $k_1, k_2, ..., k_r $ by the generating series

\begin{equation}\label{e1}
\begin{array}{c}
\frac{2Li_{(k_1,...,k_r)}(1-e^{-t})}{(1+e^t)^r}e^{rxt}=\sum\limits_{n=0}^{\infty}\mathbf{E}_{n}^{(k_1,...,k_r)}(x){t^{n}\over n!}
\end{array}.
\end{equation}
\end{defin}

and if $x=0$, then  we can define Multi Poly-Euler numbers $\mathbf{E}_{n}^{(k_1,...,k_r)}=\mathbf{E}_{n}^{(k_1,...,k_r)}(0)$

Now we define three parameters $a,b,c$, for Multi Poly-Euler polynomials and Multi Poly-Euler numbers as follows.

\begin{defin}\label{d1} Multi Poly-Euler polynomials $\mathbf{E}_n^{(k_1,...,k_r)}(x,a,b)$, $(n = 0, 1, 2, ...)$ are defined for each integer $k_1, k_2, ..., k_r $ by the generating series

\begin{equation}\label{e1}
\begin{array}{c}
\frac{2Li_{(k_1,...,k_r)}(1-(ab)^{-t})}{(a^{-t}+b^t)^r}e^{rxt}=\sum\limits_{n=0}^{\infty}\mathbf{E}_{n}^{{(k_1,...,k_r)}}(x,a,b){t^{n}\over n!}
\end{array}.
\end{equation}
\end{defin}
In the same way, and if $x=0$, then  we can define Multi Poly-Euler numbers with $a,b$ parameters $\mathbf{E}_{n}^{(k_1,...,k_r)}(a,b)=\mathbf{E}_{n}^{(k_1,...,k_r)}(0,a,b)$.

In the following theorem, we find a relation between $\mathbf{E}_{n}^{(k_1,...,k_r)}(a,b)$ and $\mathbf{E}_{n}^{(k_1,...,k_r)}(x)$

\begin{theo}\label{t1}
Let $a,b>0$, $ab\neq \pm1$ then we have

\begin{equation}\label{e1}
\begin{array}{c}
\mathbf{E}_{n}^{(k_1,k_2,...,k_r)}(a,b)=\mathbf{E}_{n}^{(k_1,k_2,...,k_r)}\left(\frac{lna}{lna+lnb}\right)(lna+lnb)^n
\end{array}.
\end{equation}

\end{theo}
Proof.By applying the Definition 2 and  Definition 3,we have

\begin{align*}
\frac{2Li_{(k_1,...,k_r)}(1-(ab)^{-t})}{(a^{-t}+b^t)^r} &= \sum\limits_{n=0}^{\infty}\mathbf{E}_{n}^{{(k_1,...,k_r)}}(a,b){t^{n}\over n!} \\
& = e^{rt\ln a}\frac{2Li_{(k_1,...,k_r)}(1-e^{-t\ln ab})}{(1+e^{t\ln ab})^r}\\
\end{align*}
So, we get
\begin{align*}
\frac{2Li_{(k_1,...,k_r)}(1-(ab)^{-t})}{(a^{-t}+b^t)^r}=\sum\limits_{n=0}^{\infty}\mathbf{E}_{n}^{{(k_1,...,k_r)}}\left(\frac{\ln a}{\ln a +\ln b}\right)(\ln a +\ln b)^n{t^{n}\over n!} 
\end{align*}
Therefore, by comparing the coefficients of $t^n$ on both sides, we get the desired result. $\square$

Now, In next theorem, we show a shortest relationship between $\mathbf{E}_{n}^{(k_1,k_2,...,k_r)}(a,b)$
and $\mathbf{E}_{n}^{(k_1,k_2,...,k_r)}$.

\begin{theo}\label{t1}
Let $a,b>0$, $ab\neq \pm1$ then we have

\begin{equation}\label{e1}
\begin{array}{c}
\mathbf{E}_{n}^{(k_1,k_2,...,k_r)}(a,b)=\sum\limits_{i=0}^{n}r^{n-i}(\ln a +\ln b)^i(\ln a)^{n-i}\binom{n}{i} \mathbf{E}_{i}^{(k_1,k_2,...,k_r)}
\end{array}.
\end{equation}

\end{theo}
Proof. By applying the Definition 2, we have,
\begin{align*}
\sum\limits_{n=0}^{\infty}\mathbf{E}_{n}^{{(k_1,...,k_r)}}(a,b){t^{n}\over n!}&=\frac{2Li_{(k_1,...,k_r)}(1-(ab)^{-t})}{(a^{-t}+b^t)^r}   \\
& = e^{rt\ln a}\frac{2Li_{(k_1,...,k_r)}(1-e^{-t\ln ab})}{(1+e^{t\ln ab})^r}\\
& = \left(\sum\limits_{k=0}^{\infty}\frac{r^kt^k(\ln a)^k}{k!}\right)\left(\sum\limits_{n=0}^{\infty}\mathbf{E}_{n}^{{(k_1,...,k_r)}}(\ln a+\ln b)^n\frac{t^n}{n!}\right)\\
& = \sum\limits_{j=0}^{\infty}\left( \sum\limits_{i=0}^{j}r^{j-i}\frac{\mathbf{E}_j^{(k_1,...,k_r)}(\ln a+\ln b)^i(\ln a)^{j-i}}{i!(j-i)!}t^j\right)\\
\end{align*}
So, by comparing the coefficients of $t^n$ on both sides , we get the desired result. $\square$

By applying the definition 2, by simple manipulation, we get the following corollary

\begin{cor}\label{c1}For non-zero numbers $a,b$, with $ab \neq -1$ we have 

\begin{equation}\label{e2}
\begin{array}{c}
\mathbf{E}_n^{(k_1,...,k_r)}(x;a,b)=\sum\limits_{i=0}^{n}\binom{n}{i}r^{n-i}\mathbf{E}_i^{(k_1,...,k_r)}(a,b)x^{n-i}
\end{array}.
\end{equation}
\end{cor}

Furthermore, by combinig the results of Theorem 2, and Corollary 1, we get the following relation between generalization of Multi Poly-Euler polynomials with $a,b$ parameters $\mathbf{E}_{n}^{{(k_1,...,k_r)}}(x;a,b)$, and  Multi Poly-Euler numbers $\mathbf{E}_{n}^{{(k_1,...,k_r)}}$. 

\begin{equation}\label{e1}
\begin{array}{c}
\mathbf{E}_{n}^{{(k_1,...,k_r)}}(x;a,b)= \sum\limits_{k=0}^{n} \sum\limits_{j=0}^{k}r^{n-k}\binom{n}{k}\binom{k}{j}(\ln a)^{k-j}(\ln a+\ln b)^j\mathbf{E}_{j}^{{(k_1,...,k_r)}}x^{n-k}
\end{array}.
\end{equation}

Now, we state the ''Addition formula'' for generalized Multi Poly-Euler polynomials

\begin{cor}\label{c1}(Addition formula) For non-zero numbers $a,b$, with $ab \neq -1$ we have 

\begin{equation}\label{e2}
\begin{array}{c}
\mathbf{E}_n^{(k_1,...,k_r)}(x+y;a,b)=\sum\limits_{k=0}^{n}\binom{n}{k}r^{n-k}\mathbf{E}_k^{(k_1,...,k_r)}(x;a,b)y^{n-k}
\end{array}.
\end{equation}
\end{cor}
Proof. We can write 

\begin{align*}
\sum\limits_{n=0}^{\infty}\mathbf{E}_{n}^{{(k_1,...,k_r)}}(x+y;a,b){t^{n}\over n!}&=\frac{2Li_{(k_1,...,k_r)}(1-(ab)^{-t})}{(a^{-t}+b^t)^r} e^{(x+y)rt} \\
& = \frac{2Li_{(k_1,...,k_r)}(1-(ab)^{-t})}{(a^{-t}+b^t)^r} e^{xrt}e^{yrt}\\
& = \left(\sum\limits_{n=0}^{\infty}\mathbf{E}_{n}^{{(k_1,...,k_r)}}(x;a,b){t^{n}\over n!}\right)
\left(\sum\limits_{i=0}^{n}\frac{y^i r^i}{i!}t^i\right)
\\
& = \sum\limits_{n=0}^{\infty}\left( \sum\limits_{k=0}^{n}\binom{n}{k}r^{n-k}y^{n-k}\mathbf{E}_{k}^{{(k_1,...,k_r)}}(x;a,b)\right)\frac{t^n}{n!}\\
\end{align*}
So, by comparing the coefficients of $t^n$ on both sides , we get the desired result. $\square$

\section{Explicit formula for Multi Poly-Euler polynomials}

Here we present an explicit formula for Multi Poly-Euler polynomials.

\begin{theo}\label{t1}
The Multi Poly-Euler polynomials have the following explicit formula
\begin{equation}\label{e2}
\begin{array}{c}
\mathbf{E}^{(k_1,k_2,\ldots, k_r)}_n(x)=\sum\limits_{i=0}^n\sum\limits_{ 0\le m_1\le m_2\le\ldots\le m_r \atop c_1+c_2+\ldots=r}\sum\limits_{j=0}^{m_r}\frac{2(rx-j)^{n-i}r!(-1)^{j+c_1+2c_2+\ldots}(c_1+2c_2+\ldots)^i\binom{m_r}{j}\binom{n}{i}}{(c_1!c_2!\ldots)(m_1^{k_1} m_2^{k_2}\ldots m_r^{k_r})}.
\end{array}.
\end{equation}
\end{theo}
Proof. We have
$$Li_{(k_1,k_2,\ldots, k_r)}(1-e^{-t})e^{rxt}=\sum_{ 0\le m_1\le m_2\le\ldots\le m_r }\frac{(1-e^{-t})^{m_r}}{m_1^{k_1} m_2^{k_2}\ldots m_r^{k_r}}e^{rxt}\qquad\qquad\qquad\qquad\qquad\qquad$$
\begin{align*}
=&\sum_{ 0\le m_1\le m_2\le\ldots\le m_r }\frac{1}{m_1^{k_1} m_2^{k_2}\ldots m_r^{k_r}}\sum_{j=0}^{m_r}(-1)^j\binom{m_r}{j}\sum_{n\ge0}(rx-j)^n\frac{t^n}{n!}\\
=&\sum_{n\ge0}\left(\sum_{ 0\le m_1\le m_2\le\ldots\le m_r }\sum_{j=0}^{m_r}\frac{(-1)^{j}(rx-j)^n\binom{m_r}{j}}{m_1^{k_1} m_2^{k_2}\ldots m_r^{k_r}}\right)\frac{t^n}{n!}.
\end{align*}
On the other hand, 
\begin{align*}
\left(\frac{1}{1+e^{t}}\right)^r=&\left(\sum_{ n\ge0 }(-1)^ne^{nt}\right)^r\\
=&\sum_{c_1+c_2+\ldots=r}\frac{r!(-1)^{c_1+2c_2+\ldots}}{c_1!c_2!\ldots}e^{t(c_1+2c_2+\ldots)}\\
=&\sum_{c_1+c_2+\ldots=r}\frac{r!(-1)^{c_1+2c_2+\ldots}}{c_1!c_2!\ldots}\sum_{n\ge0}(c_1+2c_2+\ldots)^n\frac{t^n}{n!}\\
=&\sum_{n\ge0}\left(\sum_{c_1+c_2+\ldots=r}\frac{r!(-1)^{c_1+2c_2+\ldots}(c_1+2c_2+\ldots)^n}{c_1!c_2!\ldots}\right)\frac{t^n}{n!}.
\end{align*}
Hence,
$$\frac{2Li_{(k_1,k_2,\ldots, k_r)}(1-e^{-t})}{(1+e^{t})^r}e^{rxt}=2Li_{(k_1,k_2,\ldots, k_r)}(1-e^{-t})e^{rxt}\left(\frac{1}{1+e^{t}}\right)^r\qquad\qquad\qquad\qquad\qquad\qquad$$
\begin{align*}
=\left(\sum_{n\ge0}\left(\sum_{ 0\le m_1\le m_2\le\ldots\le m_r }\sum_{j=0}^{m_r}\frac{(-1)^{j}(rx-j)^n\binom{m_r}{j}}{m_1^{k_1} m_2^{k_2}\ldots m_r^{k_r}}\right)\frac{t^n}{n!}\right)\times\\
\;\;\;\;\times\left(\sum_{n\ge0}\left(\sum_{c_1+c_2+\ldots=r}\frac{r!(-1)^{c_1+2c_2+\ldots}(c_1+2c_2+\ldots)^n}{c_1!c_2!\ldots}\right)\frac{t^n}{n!}\right)\\
=2\sum_{n\ge0}\sum_{i=0}^n\left(\sum_{ 0\le m_1\le m_2\le\ldots\le m_r }\sum_{j=0}^{m_r}\frac{(-1)^{j}(rx-j)^{n-i}\binom{m_r}{j}}{m_1^{k_1} m_2^{k_2}\ldots m_r^{k_r}}\right)\frac{t^{n-i}}{(n-i)!}\times\\
\;\;\;\;\times\left(\sum_{c_1+c_2+\ldots=r}\frac{r!(-1)^{c_1+2c_2+\ldots}(c_1+2c_2+\ldots)^i}{c_1!c_2!\ldots}\right)\frac{t^i}{i!}\\\end{align*}
\begin{align*}=2\sum_{n\ge0}\sum_{i=0}^n\sum_{ 0\le m_1\le m_2\le\ldots\le m_r \atop c_1+c_2+\ldots=r}\sum_{j=0}^{m_r}\frac{(rx-j)^{n-i}r!(-1)^{j+c_1+2c_2+\ldots}(c_1+2c_2+\ldots)^i\binom{m_r}{j}\binom{n}{i}}{(c_1!c_2!\ldots)(m_1^{k_1} m_2^{k_2}\ldots m_r^{k_r})}\frac{t^{n}}{n!}\end{align*}

By comparing the coefficient of $t^n/n!$, we obtain the desired explicit formula.

\begin{defin}\label{d1}
(Poly-Euler polynomials with $a,b,c$ parameters):The Poly-Euler polynomials with $a,b,c$ parameters may be defined by using the following generating function,

\begin{equation}\label{e1}
\begin{array}{c}
\frac{2Li_k(1-(ab)^{-t})}{a^{-t}+b^t}c^{xt}=\sum\limits_{n=0}^{\infty}\mathbf{E}_{n}^{(k)}(x;a,b,c){t^{n}\over n!}
\end{array}.
\end{equation}
\end{defin}

Now, in next theorem, we give an explicit formula for Poly-Euler polynomials with $a,b,c$ parameters.

\begin{theo}\label{t1}
The generalized Poly-Euler polynomials with $a,b,c$ parameters have the following
explicit formula

\begin{equation}\label{e2}
\begin{array}{c}
\mathbf{E}^{(k)}_n(x;a,b,c)=\\
\sum\limits_{m=0}^n\sum\limits_{j=0}^m\sum\limits_{i=0}^j\frac{2(-1)^{m-j+i}}{j^k}\binom{j}{i}(x\ln c-(m-j+i+1)\ln a-(m-j+i+1)\ln b)^n.
\end{array}.
\end{equation}
\end{theo}
Proof. We can write
\begin{align*}\sum_{n\ge0}\mathbf{E}^{(k)}_n(x;a,b,c)\frac{t^n}{n!}=\frac{2Li_k(1-(ab)^{-t})}{a^{-t}((ab)^{-t}+1)}c^{xt}
=2a^{-t}\left(\sum_{n\ge0}(-1)^n(ab)^{-nt}\right)\left(\sum_{n\ge0}\frac{\left(1-(ab)^{-t}\right)^m}{m^k}\right)c^{xt}.\end{align*}
\begin{align*}
=&a^{-t}\sum_{m\ge0}\sum_{j=0}^m\sum_{i=0}^j\frac{2(-1)^{m-j+i}}{j^k}\binom{j}{i}(ab)^{-t(x+m-j+i)}c^{xt}\\
=&\sum_{m\ge0}\sum_{j=0}^m\sum_{i=0}^j\frac{2(-1)^{m-j+i}}{j^k}\binom{j}{i}e^{-t(x+m-j+i)\ln(ab)}e^{-t\ln a}e^{xt\ln c}\\
=&\sum_{n\ge0}\sum_{m\ge0}\sum_{j=0}^m\sum_{i=0}^j\frac{2(-1)^{m-j+i}}{j^k}\binom{j}{i}\sum_{n\ge0}(x\ln c-(m-j+i+1)\ln a-(m-j+i)\ln b)^n\frac{t^n}{n!}\\
=&\sum_{n\ge0}\sum_{m=0}^n\sum_{j=0}^m\sum_{i=0}^j\frac{2(-1)^{m-j+i}}{j^k}\binom{j}{i}(x\ln c-(m-j+i+1)\ln a-(m-j+i)\ln b)^n\frac{t^n}{n!}.
\end{align*}

By comparing the coefficient of $t^n/n!$, we obtain the desired explicit formula.$\square$

\bigskip
\noindent
 $\begin{array}{ll}
\textrm{\bf Hassan Jolany}\\
\textrm{Université des Sciences et Technologies de Lille}\\
\textrm{UFR de Mathématiques}\\
\textrm{Laboratoire
Paul Painlevé}\\
\textrm{CNRS-UMR 8524
59655 Villeneuve d'Ascq Cedex/France}\\
\textrm{e-mail: hassan.jolany@math.univ-lille1.fr}
\end{array}$

\bigskip
\noindent
 $\begin{array}{ll}
\textrm{\bf Mohsen Aliabadi}\\
\textrm{Department of Mathematics, Statistics and Computer Science,}\\
\textrm{University of Illinois at Chicago, USA}\\
\textrm{e-mail: mohsenmath88@gmail.com}
\end{array}$

\bigskip
\noindent
 $\begin{array}{ll}
\textrm{\bf Roberto B. Corcino}\\
\textrm{Department of Mathematics}\\
\textrm{Mindanao State University, Marawi City, 9700 Philippines}\\
\textrm{e-mail: rcorcino@yahoo.com}
\end{array}$

\bigskip
\noindent
 $\begin{array}{ll}
\textrm{\bf M.R.Darafsheh}\\
\textrm{Department of Mathematics, Statistics and Computer Science }\\
\textrm{Faculty of Science}\\
\textrm{University of Tehran, Iran}\\
\textrm{e-mail: darafsheh@ut.ac.ir}
\end{array}$
\end{document}